\documentclass[11pt]{amsart}
\usepackage{amsfonts}
\usepackage{amsmath}
\usepackage{amssymb}
\usepackage{amsthm}
\usepackage{color,ulem}
\usepackage{mathtools}
\usepackage{enumerate}
\usepackage[a4paper,margin=1.15in]{geometry}
\usepackage[dvipsnames]{xcolor}
\usepackage{soul}
\usepackage[backref=page]{hyperref}
\usepackage{comment}
\usepackage{tikz-cd}

\input xy
\xyoption {all}

\newcommand{\bZ}{\mathbb{Z}}
\newcommand{\bP}{\mathbb{P}}
\newcommand{\bC}{\mathbb{C}}

\newcommand{\bA}{\mathbb{A}}

\newcommand{\bQ}{{\mathbb Q}}

\newcommand{\bH}{{\mathbb H}}

\def\L{\Lambda}

\def\O{\Omega}
 
\def\cC{{\mathcal C}}

\def\cE{{\mathcal E}}
\def\cF{{\mathcal F}}
\def\cG{{\mathcal G}}

\def\cI{{\mathcal I}}

\def\cK{{\mathcal K}}

\def\cM{{\mathcal M}}

\def\cO{{\mathcal O}}
\def\cP{{\mathcal P}}

\def \Ext {\operatorname{Ext}} 
 
\def \ext {{\mathcal {E}xt}}

\def\codim{\operatorname{codim}}

\def\Num{\operatorname{Num}}
\def\Pic{\operatorname{Pic}}
\newcommand{\IH}{I\!H}
\newcommand{\RIH}{RI\!H}
\newcommand{\IC}{\operatorname{IC}}

\theoremstyle{plain}
\newtheorem{thm}{Theorem}[section]
\newtheorem{cor}[thm]{Corollary}
\newtheorem{lem}[thm]{Lemma}
\newtheorem{prop}[thm]{Proposition}
\newtheorem{conj}[thm]{Conjecture}
\newtheorem{question}[thm]{Question}

\theoremstyle{definition}
\newtheorem{defn}[thm]{Definition}
\newtheorem{example}[thm]{Example}

\theoremstyle{definition}
\newtheorem{rmk}[thm]{Remark}


\usepackage{microtype}
\begin{document}
\baselineskip=16.25pt

\title[Stabilization of intersection Betti numbers]{Stabilization of intersection Betti numbers for moduli spaces of one-dimensional sheaves on surfaces}

\author{Fei Si}
\address{The School of Mathematics and Statistics, Xi’an Jiaotong University, 28 West Xianning Road, Xi’an, Shaanxi, P.R.China 710049}
\email{sifei@xjtu.edu.cn}

\author{Feinuo Zhang}
\address{Shanghai Center for Mathematical Sciences, Fudan University, Jiangwan Campus, Shanghai, 200438, China}
\email{fnzhang17@fudan.edu.cn}

\begin{abstract}{In this paper, we develop a unified approach to study the intersection Betti numbers of moduli spaces of one-dimensional semistable sheaves on smooth projective surfaces. Assuming the irreducibility of such moduli spaces, we prove that their intersection Betti numbers in a certain range of degrees coincide with the stable Betti numbers of Hilbert schemes of points. As an application, for surfaces with nef anticanonical divisor, we show that these intersection Betti numbers stabilize in each fixed degree, which fits into the broader context of stable cohomology for moduli spaces of sheaves; if in addition the moduli spaces are smooth, we also prove a refined stabilization result on perverse Hodge numbers.}
\end{abstract}

\maketitle

\tableofcontents

\section{Introduction}
\label{sec:intro}

In recent years, moduli spaces of one-dimensional sheaves on surfaces have attracted considerable attention. On the one hand, their connections with enumerative geometry and the theory of Higgs bundles have led to a number of new conjectures and results (e.g.,~\cite{Bou22, Bou23, GW25, KLMP24, KPS23, MS23, MT18, Obe24, PSSZ24, Yua23}). On the other hand, these moduli spaces exhibit not only asymptotic behavior analogous to that of moduli spaces of torsion-free sheaves, but also pathologies in contrast to the torsion-free case (see \cite[\S 5]{SZ25}), revealing rich geometry that remains to be understood.

Let $S$ be a connected, smooth, complex projective surface with a fixed ample divisor $H$. For a nonzero effective divisor $\beta$ on $S$ and an integer $\chi$, let $M_{\beta,\chi}$ be the coarse moduli space parametrizing S-equivalence classes of one-dimensional semistable (with respect to $H$) sheaves $\cF$ on $S$ with fixed determinant $\det(\cF)\cong\cO_S(\beta)$ and Euler characteristic $\chi(\cF)=\chi$. 
The moduli space $M_{\beta,\chi}$ admits a Hilbert--Chow morphism 
\[h:M_{\beta,\chi}\to |\beta|\]
to the linear system $|\beta|$, defined by taking the Fitting support.
Over the open subscheme $|\beta|^{\mathrm{int}}\subset |\beta|$ parametrizing integral curves, the fibers of $h$ are compactified Jacobians. 

A natural step toward understanding the geometry of $M_{\beta,\chi}$ is to determine its topological invariants. Since $M_{\beta,\chi}$ is singular in general, we consider its intersection cohomology, which is better-behaved than the ordinary cohomology. In this paper, we study the intersection Betti numbers of $M_{\beta,\chi}$, defined as the dimensions of the intersection cohomology groups $\IH^k(M_{\beta,\chi})$ (with $\bQ$-coefficients) for $k\in\bZ_{\ge0}$. Our first main result expresses these intersection Betti numbers, in a certain range of degrees, as the coefficients $b_k^\infty$ of the formal power series 
\begin{equation}
\label{eq:prodform}
   \sum_{k\ge0} b_k^\infty q^k=\prod_{m\ge1}\frac{(1+q^{2m-1})^{b_1(S)}(1+q^{2m+1})^{b_1(S)}}{(1-q^{2m})^{b_2(S)+1}(1-q^{2m+2})},
\end{equation}
where $b_i(S)$ denotes the $i$-th Betti number of $S$.

\begin{thm}
\label{thm:main}
Suppose that $M_{\beta,\chi}$ is irreducible, and denote by 
\[c\vcentcolon=\codim(M_{\beta,\chi} \setminus h^{-1}(|\beta|^{\mathrm{int}}),M_{\beta,\chi})\] 
the codimension of $M_{\beta,\chi} \setminus h^{-1}(|\beta|^{\mathrm{int}})$ in $M_{\beta,\chi}$.
If 
\begin{enumerate}
    \item $\beta$ is $\max\{1,2k-2\}$-very ample (see Definition~\ref{def:kveryamp}), and
    \item $k\le\min\left\{c-1,\frac{2}{3}\dim|\beta|\right\}$,
\end{enumerate}
then 
\[\dim_\bQ\IH^k(M_{\beta,\chi})=b_k^\infty.\]
\end{thm}


By Göttsche's formula \cite{Got90}, the number $b_k^\infty$ is exactly the $k$-th Betti number $b_k(S^{[n]})$ of the Hilbert scheme $S^{[n]}$ of $n$ points on $S$ for any integer $n\ge k$. The connection between (intersection, virtual, or ordinary) Betti numbers of moduli spaces of sheaves on $S$ and $b_k^\infty$ has been observed in various settings. 
For some moduli spaces of torsion-free sheaves, this connection has been established, for instance, in \cite{CW22,Li97}. Similar results have been obtained for moduli spaces of one-dimensional sheaves on specific surfaces such as K3 surfaces \cite{MM24}, the projective plane \cite{Yua23}, and del Pezzo surfaces \cite{PSSZ24}. The above Theorem~\ref{thm:main} provides a unifying perspective on the intersection Betti numbers of moduli spaces of one-dimensional sheaves on surfaces. 

The proof of Theorem~\ref{thm:main} has two main parts. The first part, which is the key new ingredient, shows the existence of a sufficiently large smooth open subscheme of $M_{\beta,\chi}$ (see Theorem~\ref{thm:gsm} below), and the second part computes the cohomology of such a smooth open subscheme by generalizing the approach in our previous joint work with Pi and Shen \cite{PSSZ24} (see Corollary~\ref{cor}). The argument in \cite{PSSZ24} relies crucially on the smoothness of the open subscheme $h^{-1}(|\beta|^{\mathrm{int}})$ of $M_{\beta,\chi}$ when $S$ is a del Pezzo surface. 
However, $h^{-1}(|\beta|^{\mathrm{int}})$ can be singular in general, and the standard dimension estimates as in \cite[Ch.~9]{HL10} for singular loci of moduli spaces of torsion-free sheaves do not apply to $M_{\beta,\chi}$.
We overcome this by using the deformation theory of curve singularities to show that the dimension of the singular locus of $h^{-1}(|\beta|^{\mathrm{int}})$ is governed by the positivity of $\beta$, which may be of independent interest. The precise statement is as follows.

\begin{thm}
\label{thm:gsm}
    If $\beta$ is $\max\{1,2k-2\}$-very ample, then there exists an open subscheme $U_k\subset|\beta|^{\mathrm{int}}$ such that
    \begin{enumerate}
    \item $\codim(|\beta|^{\mathrm{int}}\setminus U_k,|\beta|^{\mathrm{int}})\ge k+1$;
        \item the open subscheme $h^{-1}(U_k)$ of $M_{\beta,\chi}$ is smooth.
    \end{enumerate}
\end{thm}

Motivated by Theorem~\ref{thm:main} and the results for torsion-free sheaves, we propose the following conjecture on the stabilization of intersection Betti numbers, which can be viewed as a rank $0$ analogue of \cite[Conjecture~1.1]{CW22}.
\begin{conj}
\label{conj}
Given any ample divisor $\beta_0$ on $S$ and any $k\in\bZ_{\ge 0}$, $\chi\in\bZ$, there exists $d(\beta_0,k,\chi)\in\bZ_{>0}$ such that for every integer $d\ge d(\beta_0,k,\chi)$, the intersection Betti number $\dim_\bQ \IH^k(M_{d\beta_0,\chi})$ is independent of $d$. If furthermore $M_{d\beta_0,\chi}$ is irreducible, then the stable intersection Betti number is given by
\begin{equation}
\label{eq:main}
    \dim_\bQ \IH^k(M_{d\beta_0,\chi})=b_k^\infty.
\end{equation}
\end{conj}

As an application of Theorem~\ref{thm:main}, we prove that Conjecture~\ref{conj} holds if the anticanonical divisor $-K_S$ of $S$ is nef. Surfaces with nef anticanonical divisor include minimal surfaces of Kodaira dimension $0$ (i.e., K3, Abelian, Enriques, and bielliptic surfaces) and weak del Pezzo surfaces; see \cite[Propositions~1.5, 1.6]{BP04} for more types. Note that Example ~\ref{ex:gentype} also provides some evidence for Conjecture~\ref{conj} for certain surfaces of general type.
\begin{thm}
\label{thm:stab}
    Suppose that $-K_S$ is nef. Then Conjecture~\ref{conj} is true. Moreover, $d(\beta_0,k,\chi)$ can be chosen independently of $\chi$ such that \eqref{eq:main} holds for all integers $d\ge d(\beta_0,k)\vcentcolon=d(\beta_0,k,\chi)$. 
\end{thm}

The perverse Hodge numbers $n_{\beta,\chi}^{i,j}$ for $i,j\in\bZ_{\ge0}$ are defined by \eqref{eq:nij} with $\varphi$ taken to be $h:M_{\beta,\chi}\to|\beta|$, and they refine the intersection Betti numbers. Let $n_\infty^{i,j}$ denote the coefficients of the following formal power series
    \begin{equation}
    \label{eq:nijprod}
        \sum_{i,j\ge0}n_\infty^{i,j}q^it^j=
        (1-qt)\prod_{m\ge1}\frac{(1+q^mt^{m-1})^{b_1(S)}(1+q^mt^{m+1})^{b_1(S)}}{(1-q^{m+1}t^{m-1})(1-q^mt^m)^{b_2(S)}(1-q^{m-1}t^{m+1})}.
    \end{equation}
We also prove, as a refined version of Theorem~\ref{thm:stab}, the stabilization of $n_{\beta,\chi}^{i,j}$ in terms of \eqref{eq:nijprod} when $M_{\beta,\chi}$ is smooth. 

\begin{thm}
\label{thm:Enrs}
Suppose that $-K_S$ is nef. Given any ample divisor $\beta_0$ on $S$ and any $i,j\in\bZ_{\ge0}$, there exists $d(\beta_0,i,j)\in\bZ_{>0}$ such that for every integer $d\ge d(\beta_0,i,j)$, if $M_{d\beta_0,\chi}$ is smooth, then
    \begin{equation*}
    n_{d\beta_0,\chi}^{i,j}=n_\infty^{i,j}.
    \end{equation*}
\end{thm}

    When $S$ is a K3 or an Abelian surface, the statement of Conjecture~\ref{conj} for $\chi=1$ can be deduced from known results in the literature, such as those in \cite{Yos01}, and Theorem~\ref{thm:Enrs} can be deduced from \cite{dCMS22,SY22}. In the K3 case, Mauri and Migliorini \cite{MM24} proved the stabilization of intersection Betti and intersection Hodge numbers of $M_{\beta,\chi}$ for $\chi\neq 1$; the Abelian case can be proved similarly using the argument in \cite[\S 10.3]{MM24}. As in \cite{MM24,PSSZ24}, an analysis of $\codim(|d\beta_0|\setminus|d\beta_0|^{\mathrm{int}},|d\beta_0|)$ is needed for our proofs of Theorems~\ref{thm:stab} and \ref{thm:Enrs}. However, our approach to determining the stable numbers is different from that in \cite{MM24} and does not involve any hyper-Kähler geometric properties as used in \cite{dCMS22,MM24,SY22,Yos01}. Indeed, the deformation equivalence of moduli spaces with Hilbert schemes of points on $S$ needed in \cite[\S 10.3]{MM24} fails in the Enriques case (since in this case  $M_{\beta,\chi}$ is odd-dimensional, e.g., \cite{Sac19}) and is unknown in other cases. The del Pezzo case for Theorems~\ref{thm:stab} and \ref{thm:Enrs} was proved in \cite{PSSZ24} in a slightly different form.

Although moduli spaces of sheaves on K3 and Abelian surfaces have been extensively studied, much less is known about the cases of other surfaces with nef anticanonical divisor. In the bielliptic case, some progress on the birational geometry and singularities of moduli spaces was made in \cite{Nue25}. The Hodge numbers of some moduli spaces of stable one-dimensional sheaves on an elliptic ruled surface were computed in \cite{Yos23}.

When $S$ is a generic Enriques surface, assuming that $\beta^2>0$ and the class of $\beta$ is not divisible by $2$ in the group $\Num(S)$ of numerical equivalence classes of divisors on $S$, Saccà \cite{Sac19} proved that $M_{\beta,1}$ is a smooth projective Calabi--Yau variety of dimension $\beta^2+1$ and computed its first few Betti numbers. In this case, although the odd cohomology of $M_{\beta,1}$ can be nonzero as shown in \cite[Theorem~5.10]{Sac19}, our Theorem~\ref{thm:stab} implies that for fixed odd $k$, the $k$-th Betti number $b_k(M_{\beta,1})$ vanishes if $\beta$ is sufficiently positive. The invariants $n_{\beta,1}^{i,j}$ are conjectured to refine the Gromov--Witten/Pandharipande--Thomas invariants of the local surface $\operatorname{Tot}(K_S)$ \cite{HST01,KL12,MT18} and are the shifted perverse Hodge numbers considered in \cite{Obe24}.

\begin{rmk}\label{rmk:Obe}
Under the assumption of \cite[Conjecture~A]{Obe24}, Oberdieck proved a stronger stabilization result \cite[Proposition~1.1]{Obe24} for $n_{\beta,1}^{i,j}$ in the case of generic Enriques surfaces. We note that in this case, the product formula \eqref{eq:nijprod} coincides with the one in \cite[Proposition~1.1]{Obe24}. It remains open whether our approach, which does not rely on \cite[Conjecture~A]{Obe24}, can be modified to show $n_{\beta,1}^{i,j}=n^{i,j}_\infty$ if $i<\beta^2/4$ and $j<\beta^2/4-1$ as in \cite[Proposition~1.1]{Obe24}.
An explicit expression for $d(\beta_0,i,j)$ in Theorem~\ref{thm:Enrs} when $S$ is a generic Enriques surface will be given in \eqref{eq:d0}.
\end{rmk}

The rest of this paper is organized as follows. In Section~\ref{sec:pre}, we collect some notation, definitions, and results needed later. Section~\ref{sec:def} is devoted to the local versality of the universal family of curves in high codimension and contains the proof of Theorem~\ref{thm:gsm}. In Section~\ref{sec:IH}, we prove Theorem~\ref{thm:main} and discuss the geometric conditions appearing in its statement. Finally, in Section~\ref{sec:stab}, we prove Theorems~\ref{thm:stab} and \ref{thm:Enrs}, with a more explicit treatment of the Enriques case, and conclude with a discussion of related open questions.



\subsection*{Conventions}
All schemes are assumed to be of finite type over the field $\bC$ of complex numbers.
A variety is an integral and separated scheme. A curve is a projective scheme of pure dimension $1$. A surface is a variety of dimension $2$. By a point of a scheme, we always mean a closed point. All sheaves are assumed to be coherent. For a sheaf $\cE$ on a projective scheme $X$, we denote by $h^i(\cE)$ (resp. $\operatorname{ext}^i(\cE,\cE)$) the dimension of the $i$-th sheaf cohomology $H^i(X,\cE)$ (resp. $i$-th Ext group $\Ext^i(\cE,\cE)$) as a $\bC$-vector space.

\section{Preliminaries}
\label{sec:pre}

\subsection{Moduli spaces of one-dimensional sheaves}
We first review some basic notions concerning the moduli spaces studied in this paper and then discuss some properties regarding the Hilbert--Chow morphism and irreducibility. 

Fix a smooth polarized surface $(S,H)$ as in Section~\ref{sec:intro}. Recall that the dimension of a sheaf is the dimension of its support. A one-dimensional sheaf $\cF$ on $S$ is called semistable (resp. stable) (with respect to $H$) if every nonzero subsheaf $\cG\subsetneq\cF$ is one-dimensional with
\[\frac{\chi(\cG)}{c_1(\cG)\cdot H}\le(\text{resp.}<)\, \frac{\chi(\cF)}{c_1(\cF)\cdot H}.\]
Two semistable sheaves are called S-equivalent if their Jordan--Hölder filtrations have isomorphic factors up to reordering. Each S-equivalence class is represented by a unique polystable sheaf, i.e., a semistable sheaf which is the direct sum of stable sheaves.
In general, the coarse moduli space $M_{\beta,\chi}$ considered in this paper depends on the ample divisor $H$, but we have omitted $H$ from the notation since $H$ is fixed throughout.

Let $\cM_{\beta,\chi}$ be the moduli stack of one-dimensional semistable sheaves $\cF$ on $S$ with $\det(\cF)\cong\cO_S(\beta)$ and $\chi(\cF)=\chi$. Choose a large integer $n$ such that for any semistable sheaf $\cF$ in $\cM_{\beta,\chi}$, the sheaf $\cF(nH)$ is globally generated and $h^j(\cF(nH))=0$ for $j>0$. Set $N\vcentcolon=(\beta\cdot H)n+\chi$. Denote by $Q$ the Quot scheme parametrizing quotients $\cO_S(-nH)^{\oplus N}\to \cE$, where $\cE$ is semistable with fixed Hilbert polynomial $\chi(\cE(mH))=(\beta\cdot H)m+\chi$ and the induced map $H^0(S,\cO_S^{\oplus N})\to H^0(S,\cE(nH))$ is an isomorphism. 

Let $\Pic(S)$ denote the Picard scheme of $S$. There is a morphism
\[\det:Q\to\Pic(S)\]
sending $[\cO_S(-nH)^{\oplus N}\to \cE]$ to $\det(\cE)$. Then $\cM_{\beta,\chi}$ has a quotient stack presentation
\[\cM_{\beta,\chi}\simeq[Q_{\beta}/\operatorname{GL}(N)],\]
where $Q_\beta\vcentcolon=\det^{-1}(\cO_S(\beta))$ and $\operatorname{GL}(N)$ is the general linear group acting naturally on $Q_\beta$. 

The natural morphism
\[\nu:\cM_{\beta,\chi}\to M_{\beta,\chi}\]
is the good moduli space for $\cM_{\beta,\chi}$. Denote by
\[\tilde h=h\circ\nu:\cM_{\beta,\chi}\to |\beta|\] 
the stacky Hilbert--Chow morphism. Since the fibers of $h$ over the locus $|\beta|^{\mathrm{int}}$ of integral curves are compactified Jacobians and have dimension $p_a(\beta)$, where $p_a(\beta)$ is the arithmetic genus of $\beta$, the fibers of $\tilde h$ over $|\beta|^{\mathrm{int}}$ have dimension $p_a(\beta)-1$. Whether the fiber dimension of $\tilde h$ jumps depends on numerical properties of the canonical divisor $K_S$. 

\begin{lem}[{\cite[Corollary~1.3]{Yua23}}]
\label{lem:fibdim}
  If $-K_S$ is nef and $|\beta|^{\mathrm{int}}\neq\varnothing$, then for any curve $C\in|\beta|$, the fiber $\tilde h^{-1}(C)$ has dimension $p_a(\beta)-1$.
\end{lem}

This property of $\tilde h$ gives constraints on the geometry of $M_{\beta,\chi}$. The next proposition generalizes \cite[Theorem~2.3]{MS23} and \cite[Theorem~1.5]{Yua23} beyond the del Pezzo case. 
\begin{prop}
\label{prop:irr}
Suppose that $-K_S$ is nef and $|\beta|$ contains a smooth connected curve. If $\dim|\beta|=\chi(\cO_S(\beta))-1$, then
\begin{enumerate}[(i)]
    \item the moduli space $M_{\beta,\chi}$ is irreducible;
    \item all fibers of $h:M_{\beta,\chi}\to|\beta|$ have dimension $p_a(\beta)$. 
\end{enumerate}
\end{prop}
\begin{proof}
(i) It suffices to show that $Q_\beta$ is irreducible. We use the deformation theory of sheaves and quotients (see \cite[Theorem~6.4.9, Proposition~6.5.1]{FGIKNV} and \cite[2.A.8, Theorem~4.5.3]{HL10}). For any $\xi\in Q_{\beta}$ representing the short exact sequence
\begin{equation}
\label{eq:ses}
    0\to\cK\to\cO_S(-nH)^{\oplus N}\to \cF\to0,
\end{equation}
the dimension of $Q_{\beta}$ at $\xi$, denoted by $\dim_\xi Q_{\beta}$, satisfies (cf.~\cite[Proposition~4.5.9]{HL10})
\begin{equation}
\label{eq:irr}
\begin{aligned}
  \dim_\xi Q_{\beta}&\ge\dim_\xi Q-\dim\Pic(S)\\
 &\ge \operatorname{ext}^0(\cK,\cF)-(\operatorname{ext}^2(\cF,\cF)-h^2(\cO_S))-h^1(\cO_S)\\
 &=N^2+\beta^2+\chi(\cO_S)-1,
\end{aligned}
\end{equation}
where the last equality follows from the long exact sequence of Ext groups for \eqref{eq:ses} and the Riemann--Roch formula.
By \cite[Proposition~2.2]{SZ25}, the open subscheme $h^{-1}(|\beta|^{\mathrm{int}})$ is irreducible and contains the smooth open subscheme $h^{-1}(|\beta|^{\mathrm{sm}})$, where $|\beta|^{\mathrm{sm}}\subset|\beta|$ is the open subscheme parametrizing smooth curves. Then by \cite[Corollary~4.3.5]{HL10}, the open subscheme $Q_\beta^{\mathrm{int}}\subset Q_{\beta}$ parametrizing those quotient sheaves with integral Fitting support is irreducible and contains the smooth open subscheme $Q_\beta^{\mathrm{sm}}$ parametrizing those quotient sheaves with smooth Fitting support. There is a unique irreducible component $Y\subset Q_{\beta}$ containing $Q_\beta^{\mathrm{int}}$. If $Y\not=Q_{\beta}$, by taking an irreducible component $Y'\neq Y$ and noting that $Y'\cap Q_\beta^{\mathrm{sm}}=\varnothing$, we deduce from  Lemma~\ref{lem:fibdim} and \cite[Tag 0DS4]{Sta26} that 
\begin{equation*}
\dim Y'<\dim\operatorname{GL}(N)+\dim|\beta|+(p_a(\beta)-1)=N^2+\beta^2+\chi(\cO_S)-1,
\end{equation*}
which contradicts \eqref{eq:irr}. Hence, $Q_{\beta}=Y$ is irreducible, and so is $M_{\beta,\chi}$.

(ii) 
By the upper semicontinuity of the fiber dimension of $h$, $\dim h^{-1}(C)\ge p_a(\beta)$ for any curve $C\in|\beta|$. 
Note that $\cM_{\beta,\chi}$ is also irreducible since $Q_\beta$ is.
Consider the morphism $\nu:\cM_{\beta,\chi}\to M_{\beta,\chi}$ and its restriction $\nu|_{\tilde h^{-1}(C)}:\tilde h^{-1}(C)\to h^{-1}(C)$. By Lemma~\ref{lem:fibdim} and \cite[Tag 0DS4]{Sta26}, we have
\[\dim\tilde h^{-1}(C)-\dim h^{-1}(C)=(p_a(\beta)-1)-\dim h^{-1}(C) \ge -1=\dim\cM_{\beta,\chi}-\dim M_{\beta,\chi},\]
which forces $\dim h^{-1}(C)=p_a(\beta)$.
\end{proof}

\begin{rmk}
If $\beta$ is an ample divisor on a K3 (resp. an Abelian) surface $S$ satisfying $\codim(|\beta| \setminus |\beta|^{\mathrm{int}},|\beta|)\ge2$, Proposition~\ref{prop:irr} (ii) can also be obtained from the holomorphic symplectic geometry of $M_{\beta,\chi}$:\footnote{The authors are grateful to Mirko Mauri for pointing out this alternative approach and the relevant references via private communication.} Using \cite[Propositions~3.5, 3.14]{BCGPSV}, one can show that $M_{\beta,\chi}$ (resp. $K_{\beta,\chi}$) is a primitive symplectic variety (e.g., \cite[Definition~3.3]{BCGPSV}) admitting a Lagrangian fibration over $|\beta|$ (with no restrictions on $H$ or $\chi$), and thus all fibers of $h:M_{\beta,\chi}\to|\beta|$ have dimension $p_a(\beta)$ by \cite[Theorem~3]{Sch20}. Here,
when $S$ is an Abelian surface with identity $0_S$, we denote by $K_{\beta,\chi}$ the fiber $\widehat{\varphi}^{\,-1}(0_S)$ of the isotrivial fibration $\widehat{\varphi}: M_{\beta,\chi} \to S$ (cf.~\cite[Definition~4.1]{Yos01}).
\end{rmk}

\subsection{Intersection cohomology and perverse filtrations}
For a quasi-projective variety $X$, let $\IC_X$ denote the intersection cohomology complex of $X$ with $\bQ$-coefficients and middle perversity. The intersection cohomology groups of $X$ are defined as the hypercohomology groups of $\IC_X$: 
\[\IH^k(X)\vcentcolon=\bH^{k-\dim X}(X,\IC_X).\]
Note that $\IH^k(X)$ can also be defined if $X$ is allowed to be non-reduced (by the topological invariance of intersection cohomology) or reducible (even not pure dimensional, as shown in \cite[\S 4.6]{dCa12}).
When $X$ is smooth, $\IH^k(X)$ can be identified with the rational cohomology group $H^k(X)$.

Let $\varphi: X \to Y$ be a proper morphism between quasi-projective varieties. We assume that $\dim X = a$, $\dim Y = b$, and all fibers of $\varphi$ have dimension $a-b$. The perverse filtration
\[
P_0\IH^k(X) \subset P_1\IH^k(X) \subset \dots \subset   \IH^k(X)
\]
 induced by $\varphi$ is an increasing filtration on $\IH^k(X)$, governed by the topology of the morphism $\varphi$. 
It is defined by
\[
P_i\IH^k(X) \vcentcolon= \mathrm{Im}\left\{ \bH^{k-b}(Y, {^\mathbf{p}\tau_{\leq i}} (R\varphi_* \IC_X[b-a])) \to \bH^{k-b}(Y, R\varphi_* \IC_X[b-a])\right\},
\]
where $^\mathbf{p}\tau_{\leq \bullet}$ is the perverse truncation functor. We denote the dimension of the graded piece $\mathrm{Gr}_i^P \IH^{i+j}(X)\vcentcolon=P_i \IH^{i+j}(X)/P_{i-1} \IH^{i+j}(X)$ by  
\begin{equation}
\label{eq:nij}
   n_{\varphi}^{i,j}\vcentcolon=\dim_\bQ \mathrm{Gr}_i^P \IH^{i+j}(X). 
\end{equation}
When $X$ is smooth, the decomposition theorem \cite{BBD81} applied to $\varphi: X\to Y$ yields
\begin{equation*}
R\varphi_*\bQ_X[b]\simeq \bigoplus_{i=0}^{2a-2b}\mathcal{P}_i[-i],
\end{equation*}
where $\mathcal{P}_i$ are semisimple perverse sheaves on $Y$; in this case, $n_{\varphi}^{i,j}=\dim_\bQ \bH^{j-b}(Y,\cP_i)$. 

For our purposes, we mainly consider the Hilbert--Chow morphism $h:M_{\beta,\chi}\to|\beta|$ and its restrictions. The invariants $n_{\varphi}^{i,j}$ with $\varphi$ taken to be (the restrictions of) $h$ will play a central role in Section~\ref{sec:IH}.

\subsection{Positivity conditions for divisors}
To study the asymptotic behavior of $M_{\beta,\chi}$, we recall the definition of $k$-very ampleness for divisors on $S$, which generalizes base-point-freeness ($0$-very ampleness) and very ampleness ($1$-very ampleness).

\begin{defn}
\label{def:kveryamp}
    An effective divisor $\beta$ on $S$ is called \textit{$k$-very ample} if for any $0$-dimensional closed subscheme $Z\subset S$ of length $k+1$, the restriction map 
    $$r_Z:H^0(S,\cO_S(\beta))\to H^0(Z,\cO_S(\beta)|_Z)$$ is surjective. We adopt the convention that every divisor is $(-1)$-very ample. 
\end{defn}

When $h^1(\cO_S(\beta))=0$, the restriction map $r_Z$ is surjective if and only if $h^1(\cI_Z\otimes\cO_S(\beta))=0$
by the long exact sequence of cohomology groups for the short exact sequence
\[0\to\cI_Z\otimes\cO_S(\beta)\to\cO_S(\beta)\to\cO_S(\beta)|_Z\to0,\]
where $\cI_Z$ is the ideal sheaf of $Z\subset S$.
Definition~\ref{def:kveryamp} immediately implies two properties, which will be used in the subsequent sections:
\begin{enumerate}
    \item If $\beta$ is $k$-very ample and $k'\le k$ $(k'\in\bZ_{\geq0})$, then $\beta$ is also $k'$-very ample.
    \item If $\beta$ is ample, then for any $k\in\bZ_{\ge0}$, there exists $d_k\in\bZ_{>0}$ such that $d\beta$ is $k$-very ample for every integer $d\ge d_k$.
\end{enumerate}

\section{Singular loci of relative compactified Jacobians}
\label{sec:def}

The goal of this section is to prove Theorem~\ref{thm:gsm}. 
We begin by reviewing some basic notions of plane curve singularities and deformation theory. 

\subsection{Singularities and deformations of locally planar curves}
\subsubsection{Singularity invariants}
We first recall some analytic invariants of isolated singularities in locally planar curves. As a main reference, we refer to \cite{GLS07}.
\begin{defn}
\label{def:inv}
    Let $f\in\bC[[x,y]]$ be a convergent power series.
    \begin{enumerate}[(i)]
    \item The numbers 
$$\mu(f)\vcentcolon=\dim_\bC\bC[[x,y]]/\left(\partial_x f,\partial_y f\right),\quad\tau(f)\vcentcolon=\dim_\bC\bC[[x,y]]/\left(f,\partial_x f,\partial_y f\right)$$
    are called the {\it Milnor number} and {\it Tjurina number} of $f$, respectively.
    \item If $f$ is reduced, let 
    $$\cO=\bC[[x,y]]/(f)\hookrightarrow\overline{\cO}$$
    be the normalization of $\cO$. Then the number
    $$\delta(f)\vcentcolon=\dim_\bC\overline{\cO}/\cO$$
    is called the \it{$\delta$-invariant} of $f$.
    \end{enumerate}
\end{defn}

The Milnor number and the $\delta$-invariant are related by the following formula of Milnor.

\begin{prop}[{\cite[\S 10]{Mil68}}]
\label{prop:Mil}
Let $f\in(x,y)\subset\bC[[x,y]]$ be a reduced, convergent power series. Denote by $r(f)$ the number of irreducible factors of $f$. Then
\begin{equation}
\label{eq:mil}
    \mu(f)=2\delta(f)-r(f)+1.
\end{equation}
\end{prop}

If $C$ is a reduced curve on $S$, then the singular locus $C_{\mathrm{sing}}$ of $C$ is discrete. For any point $p\in C_{\mathrm{sing}}$, there exists a polynomial $f\in(x,y)\subset\bC[x,y]$ without multiple factors such that the singularity of $C$ at $p$ is analytically isomorphic to that of $\{f=0\}$ at $(0,0)$, i.e.,
\[\widehat\cO_{C,p}\cong\bC[[x,y]]/(f).\]
Thus, the Milnor number, the Tjurina number, and the $\delta$-invariant of $C$ at $p$ can be defined as those of $f$ in Definition~\ref{def:inv}, denoted by $\mu_p(C)$, $\tau_p(C)$, and $\delta_p(C)$, respectively.


For completeness, we give a proof of the following standard interpretation (cf.~\cite[Ch.~II, Exercise~1.4.2]{GLS07}) of the total Tjurina number of $C$ in terms of the  Ext sheaf $\ext^1(\O_C,\cO_C)$, where $\O_C$ is the cotangent sheaf of $C$.
\begin{lem}
\label{lem:tau}
    If $C$ is a reduced curve on S, then 
    $$\sum_{p\in C_{\mathrm{sing}}}\tau_p(C)=h^0(\ext^1(\O_C,\cO_C)).$$
\end{lem}
\begin{proof}
    The sheaf $\O_C$ is locally free over the smooth locus $C\setminus C_{\mathrm{sing}}$ of $C$, so $\ext^1(\O_C,\cO_C)$ is a torsion sheaf supported on $C_{\mathrm{sing}}$. Let $p\in C_{\mathrm{sing}}$ and assume that the singularity of $C$ at $p$ is analytically isomorphic to that of $\{f=0\}$ at $(0,0)$ for $f\in(x,y)\subset\bC[x,y]$ without multiple factors. The lemma can be proved locally around $p$. For $R\vcentcolon=\bC[x,y]/(f)$, we have $\O_R\cong(R dx\oplus R dy)/(\partial_xfdx+\partial_yfdy)$. Using the conormal exact sequence
    $$0\to (f)/(f^2)\to R^{\oplus2}\to\O_{R}\to0,$$
    we obtain $\Ext^1(\O_R,R)\cong R/(\partial_xf,\partial_yf)$, from which the result follows.
\end{proof}

\subsubsection{Smoothness from local versality}
Let $\pi_B:\cC_B\to B$ be a family of reduced, locally planar curves, i.e., a proper and flat morphism whose fibers are reduced curves with at worst planar singularities. Assume that $B$ is a smooth variety. There is a natural map (cf.~ \cite[Ch.~II, Lemma~1.20, Exercise~1.4.4]{GLS07})
\begin{equation}
\label{eq:locKS}
\kappa_{\pi_B,b}^{\text{loc}}:T_bB\to \prod_{p\in (\cC_b)_{\mathrm{sing}}}T\mathrm{Def}_{\cC_b,p}=H^0(\cC_b,\ext^1(\O_{\cC_b},\cO_{\cC_b})), 
\end{equation}
called the local Kodaira--Spencer map, where $(\cC_b)_{\mathrm{sing}}$ is the singular locus of $\cC_b\vcentcolon=\pi_B^{-1}(b)$, and $T\mathrm{Def}_{\cC_b,p}$ is the tangent space of the deformation functor of the complete local $\bC$-algebra $\widehat{\cO}_{\cC_b,p}$.

\begin{defn}
\label{def:locv}
    The family $\pi_B:\cC_B\to B$ is called \textit{locally versal} if for any point $b\in B$, the local Kodaira--Spencer map $\kappa_{\pi_B,b}^{\text{loc}}$ in \eqref{eq:locKS} is surjective. 
\end{defn}

If  $\pi_B:\cC_B\to B$ is a family of integral, locally planar curves, 
then for any $e\in\bZ$, denote by
$\varrho_{B,e}:\overline{J}_B^e\to B$
the relative compactified Jacobian of degree $e$ associated with $\pi_B$. The fiber $\varrho_{B,e}^{-1}(b)$ over a point $b\in B$ is the moduli space of rank $1$, torsion-free sheaves of degree $e$ (or equivalently, with Euler characteristic $e+1-p_a(\cC_b)$) on $\cC_b$. Since $\pi_B$ admits, étale locally on $B$, a section which avoids the singularities of fibers of $\pi_B$, the smoothness of $\overline{J}_B^e$ is equivalent to that of $\overline{J}_B^0$. Hence, the following result follows directly from \cite[Corollary~B.3]{FGvS99}.

\begin{prop}[{\cite[Corollary~B.3]{FGvS99}}]
\label{prop:sm}
If $\pi_B:\cC_B\to B$ is a locally versal family of integral, locally planar curves, then $\overline{J}_B^e$ is smooth for any $e\in\bZ$.
\end{prop}

\subsection{Local versality and smoothness in high codimension}
\label{sec:smpf}
Now we show that the universal family $\pi:\cC\to |\beta|$ of curves in $|\beta|$ is locally versal in sufficiently high codimension if $\beta$ is positive enough. 

\begin{prop}
\label{prop:locv}
If $\beta$ is $\max\{1,2k-2\}$-very ample, then there exists an open subscheme $V_k\subset|\beta|$ such that
\begin{enumerate}
\item every point of $V_k$ corresponds to a reduced curve;
 \item $\codim(|\beta|\setminus V_k,|\beta|)\ge k+1$;
    \item the universal family $\pi|_{\pi^{-1}(V_k)}:\pi^{-1}(V_k)\to V_k$ of curves in $V_k$ is locally versal.
\end{enumerate}
\end{prop}
\begin{proof}
By setting $V_0\vcentcolon=V_1$, we may assume $k>0$.
The assumption on $\beta$ implies that $\beta$ is both $k$-very ample and ($2k-2$)-very ample.
We use an argument in \cite{KST11} to control the locus of curves that are not $k$-nodal. Here, a curve is called $k$-nodal if it is reduced and has no singularities other than $k$ nodes.
Since $\beta$ is $k$-very ample, by the proof of \cite[Proposition~2.1]{KST11}, there is an open subscheme $V_k\subset |\beta|$ satisfying the following conditions:
 \begin{itemize}
     \item $\codim(|\beta|\setminus V_k,|\beta|)>k$.
     \item Curves in $V_k$ are either $k$-nodal, or reduced of geometric genus $>p_a(\beta)-k$.
 \end{itemize}
Note that for a curve $C\in V_k$ which is not $k$-nodal, its geometric genus $p_g(C)$ satisfies
\begin{equation*}
p_g(C)=p_a(\beta)-\sum_{p\in C_{\mathrm{sing}}}\delta_p(C)>p_a(\beta)-k.    
\end{equation*}
From this and Milnor's formula \eqref{eq:mil}, we obtain
\begin{equation}
\label{eq:len}
    \sum_{p\in C_{\mathrm{sing}}}\tau_p(C)\le\sum_{p\in C_{\mathrm{sing}}}\mu_p(C)\le2\sum_{p\in C_{\mathrm{sing}}}\delta_p(C)<2k.
\end{equation}

It remains to show Property (3). For any curve $D\in V_k$, let $s_D\in H^0(S,\cO_S(\beta))$ be a defining section of $D$. The map $\kappa_{\pi,D}^{\text{loc}}$ as in \eqref{eq:locKS} can be identified with the map
$$\bar r_{Z_D}:H^0(S,\cO_S(\beta))/\langle s_D\rangle\to H^0(Z_D,\cO_S(\beta)|_{Z_D})$$
induced by restriction to a $0$-dimensional closed subscheme $Z_D\subset S$. By \eqref{eq:len} (which still holds for a $k$-nodal curve) and Lemma~\ref{lem:tau}, the length $h^0(\cO_{Z_D})$ of $Z_D$ satisfies
\[h^0(\cO_{Z_D})=h^0(\ext^1(\O_D,\cO_D))=\sum_{p\in D_{\mathrm{sing}}}\tau_p(D)<2k.\] Hence, $\kappa_{\pi,D}^{\text{loc}}$ is surjective by the ($2k-2$)-very ampleness of $\beta$. 
\end{proof}

Combining the above proposition with Proposition~\ref{prop:sm}, we now prove Theorem~\ref{thm:gsm}.
\begin{proof}[Proof of Theorem~\ref{thm:gsm}]
    Let $V_k$ be as in Proposition~\ref{prop:locv}. Take $U_k=V_k\cap|\beta|^{\mathrm{int}}$. By Property (2) of $V_k$ in Proposition~\ref{prop:locv}, we have
  \[\codim(|\beta|^{\mathrm{int}}\setminus U_k,|\beta|^{\mathrm{int}})\ge k+1.\]
 It remains to show that $h^{-1}(U_k)$ is smooth. Note that $h^{-1}(U_k)$ is the relative compactified Jacobian 
 of degree $p_a(\beta)-1+\chi$ associated with the universal family of curves in $U_k$. 
This universal family is locally versal by Property (3) of $V_k$ in Proposition~\ref{prop:locv}. Thus, the smoothness of $h^{-1}(U_k)$ follows from Proposition~\ref{prop:sm}.
\end{proof}

\section{Intersection Betti numbers}
\label{sec:IH}
In this section, we prove Theorem~\ref{thm:main} by combining Theorem~\ref{thm:gsm} with a generalization of the approach in \cite{PSSZ24}. To start with, we review some properties of relative Hilbert schemes of points over $|\beta|$, which will play an important role in the proof.

\subsection{Relative Hilbert schemes of points}
Let $\pi:\cC\to |\beta|$ be the universal family of curves in $|\beta|$. For $k\in\bZ_{\geq0}$, denote by $$\pi^{[k]}:\cC^{[k]}\to |\beta|$$
the relative Hilbert scheme of $k$ points on the fibers of $\pi$. For $C\in |\beta|$, the fiber of $\pi^{[k]}$ over $C$ is the Hilbert scheme $C^{[k]}$ parametrizing $0$-dimensional, length $k$ closed subschemes of $C$. Note that $\pi^{[0]}:|\beta|\to |\beta|$ is the identity and
$\pi=\pi^{[1]}:\cC\to |\beta|$.
The following property of $\pi^{[k]}$ is a corollary of \cite[Theorem~1.1]{Lua22} (see also \cite[Proposition~2.7]{SZ25}).

\begin{lem}
\label{lem:relh}
    The morphism $\pi^{[k]}:\cC^{[k]}\to|\beta|$ has fibers of the same dimension $k$.
\end{lem}

There is a natural morphism $\sigma_k:\cC^{[k+1]}\to S^{[k+1]}$ defined by
\begin{equation*}
    [Z\subset C]\mapsto [Z\subset S],
\end{equation*}
where $C\in|\beta|$, and $Z\subset C$ is a $0$-dimensional closed subscheme of length $k+1$. The fiber of $\sigma_k$ over $[Z\subset S]\in S^{[k+1]}$ is the projective space $\bP(H^0(S,\cI_Z\otimes\cO_S(\beta))^\vee)$, where $\cI_Z$ is the ideal sheaf of $Z\subset S$.
A well-known connection between $k$-very ampleness and the structure of $\cC^{[k+1]}$ as a scheme over $S^{[k+1]}$ via $\sigma_k$ is as follows (e.g., \cite[Proposition~2.6]{SZ25}).
\begin{lem}
\label{kveryamp}
    If $\beta$ is $k$-very ample, then $\cC^{[k+1]}$ is a projective bundle over $S^{[k+1]}$.
\end{lem}

Recall that for the Hilbert--Chow morphism $h:M_{\beta,\chi}\to|\beta|$, the fibers of $h$ over the locus $|\beta|^{\mathrm{int}}$ of integral curves are the compactified Jacobians of degree $p_a(\beta)-1+\chi$ and have dimension $p_a(\beta)$. For any open subscheme $U\subset |\beta|$, let $h_U:h^{-1}(U)\to U$ be the restriction of $h$, and write $\cC_U^{[k]}$ for $(\pi^{[k]})^{-1}(U)$. 
Later, we will use the following result, due to Maulik--Yun \cite[Theorem~2.13]{MY14} and Migliorini--Shende \cite[Corollary~2]{MS13}.
\begin{thm}[\cite{MY14,MS13}]
  Suppose that $|\beta|$ contains a connected smooth curve. Let $U\subset |\beta|^{\mathrm{int}}$ be an open subscheme such that $h^{-1}(U)$ is smooth, and let $P_\bullet H^*(h^{-1}(U))$ be the perverse filtration induced by $h_U:h^{-1}(U)\to U$. For any $k,\ell\in\bZ_{\ge0}$, we have an isomorphism of vector spaces\footnote{The results of \cite{MY14,MS13} apply since the smoothness of $h^{-1}(U)$ implies, via a section of the universal curve (étale locally on $U$, avoiding the singularities of fibers), that the relative compactified Jacobians $\overline{J}_{U}^e$ for all $e\in\bZ$ are smooth, and by \cite[Proposition~14]{She12}, so are the relative Hilbert schemes $\cC^{[k]}_{U}$ for all $k\in\bZ_{\geq0}$.}
    \begin{equation}
    \label{eq:supp}
       H^\ell\left(\cC^{[k]}_{U}\right)\cong\bigoplus_{i+j\le k;\,i,j\ge0}\mathrm{Gr}_i^P H^{\ell-2j}\left(h^{-1}(U)\right).
    \end{equation}
    
\end{thm}

\begin{rmk}
\label{rmk:fullsupp}
When $k=0$, the isomorphism \eqref{eq:supp} reads
    \[H^\ell(U)\cong \mathrm{Gr}_0^P H^{\ell}\left(h^{-1}(U)\right),\]
which can be obtained as a consequence of Ngô's support theorem \cite[Théorème~7.2.1]{Ngo10} (see \cite[Theorem~2.4, \S 2.5]{MY14} or \cite[Corollary~9]{MS13}). Ngô's support theorem implies that $h_U:h^{-1}(U)\to U$ has full support. More precisely, there is an isomorphism
    \[Rh_{U*}\bQ_{h^{-1}(U)}[\dim U]\simeq\bigoplus_{i=0}^{2p_a(\beta)}\operatorname{IC}(\wedge^iR^1\pi^s_*\bQ_{\cC_U^{\mathrm{sm}}})[-i],\]
    where $\pi^s:\cC_U^{\mathrm{sm}}\to U^{\mathrm{sm}}$ is the universal family of smooth curves in $U$.
\end{rmk}



\subsection{Strategy of the proof}
\label{sec:step}
For any open subscheme $U\subset|\beta|$, let $$N(U)\vcentcolon=2\codim(|\beta|\setminus U,|\beta|)-2.$$ 
We simply write $N(\beta)$ for $N(|\beta|^{\mathrm{int}})$. 
 

As a preliminary step, we focus on a smooth open subscheme $h^{-1}(U)$ of $M_{\beta,\chi}$ and prove the following result concerning graded dimensions of the perverse filtration induced by the restriction $h_U:h^{-1}(U)\to U$ of $h:M_{\beta,\chi}\to|\beta|$ (see \eqref{eq:nij}).

\begin{thm}
\label{thm}
Suppose that $|\beta|$ contains a connected smooth curve.
Let $U\subset |\beta|^{\mathrm{int}}$ be an open subscheme such that $h^{-1}(U)$ is smooth.
    For any $i,j\in\bZ_{\ge0}$ satisfying $i+j\le N(U)$ and $3i+j\le 2\dim|\beta|$, if $\beta$ is $(i-1)$-very ample, then 
    \begin{equation*}
    n_{h_U}^{i,j}=n_\infty^{i,j},
    \end{equation*}
    where $n_\infty^{i,j}$ is defined by \eqref{eq:nijprod}.
    \end{thm}
\begin{rmk}
\label{rmk}
Recall Göttsche's formula \cite{Got90} for $\sum_{k\ge0}\sum_{i\ge 0} b_i(S^{[k]})z^iw^k$:
    \begin{equation*}
G(z,w)\vcentcolon=\prod_{m\ge1} \frac{(1+z^{2m-1}w^m)^{b_1(S)}(1+z^{2m+1}w^m)^{b_1(S)}}{(1-z^{2m-2}w^m)(1-z^{2m} w^m)^{b_2(S)}(1-z^{2m+2} w^m)}.
\end{equation*}
Let $H(q,t)$ denote the right-hand side of \eqref{eq:nijprod}.
After the change of variables
    \[z=t,\quad w=\frac{q}{t},\]
    we have the following equality:
    \begin{equation}\label{eq:GH}
\frac{H(q,t)}{1-qt}=G(z,w)\cdot\frac{1-w}{1-z^2}.
\end{equation}
\end{rmk}

By \eqref{eq:nij}, the $k$-th Betti number of $h^{-1}(U)$ satisfies
\begin{equation}
\label{eq:nijsum}
    b_k(h^{-1}(U))=\sum_{i=0}^k n_{h_U}^{i,k-i}.
\end{equation}
Since the right-hand side of \eqref{eq:prodform} is exactly $H(q,q)$, 
the next corollary follows immediately from Theorem~\ref{thm} and \eqref{eq:nijsum}.

\begin{cor}
\label{cor}
In the setting of Theorem~\ref{thm}, if $k\le\min\{N(U),2\dim|\beta|/3\}$ and $\beta$ is ($k-1$)-very ample, then 
$$b_k(h^{-1}(U))=b_k^\infty,$$
where $b^{\infty}_k$ is defined by \eqref{eq:prodform}.
\end{cor}

We prove Theorem~\ref{thm} by relating $h^{-1}(U)$ to the relative Hilbert schemes $\cC^{[k]}_{U}$. By \eqref{eq:supp}, for any $k,\ell\in\bZ_{\ge0}$,
\begin{equation}
    \label{eq:myms}
       b_\ell(\cC^{[k]}_{U})=\sum_{i+j\le k;\,i,j\ge0}n_{h_U}^{i,\ell-i-2j}. 
    \end{equation}

To calculate the Betti numbers of $\cC^{[k]}_{U}$, we use Lemmas~\ref{lem:relh} and \ref{kveryamp}.

\begin{prop}
\label{prop:sk}
In the setting of Theorem~\ref{thm},
    if $\ell\le \min\{N(U),2\dim|\beta|-2k\}$ and $\beta$ is ($k-1$)-very ample, then  
\begin{equation}
\label{eq:sk}
b_\ell(\cC^{[k]}_{U})=\sum_{r=0}^{\lfloor \ell/2\rfloor}b_{\ell-2r}(S^{[k]}).   
\end{equation}
\end{prop}
\begin{proof}
Since all fibers of $\pi^{[k]}:\cC^{[k]}\to|\beta|$ have dimension $k$ by Lemma~\ref{lem:relh}, we have
$$\codim(\cC^{[k]}\setminus\cC_{U}^{[k]},\cC^{[k]})=\codim(|\beta|\setminus U,|\beta|).$$
By Lemma~\ref{kveryamp}, the $(k-1)$-very ampleness of $\beta$ implies that $\cC^{[k]}$ is a projective bundle over $S^{[k]}$. In particular, $\cC^{[k]}$ is irreducible and smooth.
    Then it follows from \cite[Lemma~2.4]{PSSZ24} (or the proof of Lemma~\ref{lem:intEnr}) that for $\ell\le N(U)$,
    $$b_\ell(\cC^{[k]}_{U})=b_\ell(\cC^{[k]}).$$
Combining this with \[\dim\cC^{[k]}-\dim S^{[k]}=\dim|\beta|-k\ge \ell/2,\] which follows from our assumption, we conclude that
$$b_\ell(\cC^{[k]}_{U})=b_\ell(\cC^{[k]})=b_\ell(S^{[k]}\times\bP^{\dim|\beta|-k})=\sum_{r=0}^{\lfloor \ell/2\rfloor}b_{\ell-2r}(S^{[k]}),$$
thereby completing the proof of \eqref{eq:sk}.
\end{proof}

Since all Betti numbers of $S^{[k]}$ are determined by Göttsche's formula, now we can prove Theorem~\ref{thm} by combining \eqref{eq:myms} and \eqref{eq:sk}.

\begin{proof}[Proof of Theorem~\ref{thm}]
We use induction on $i$. For $i=0$ and $j\le\min\{N(U),2\dim|\beta|\}$,
\[\begin{aligned}n_{h_U}^{0,j}&=\dim_\bQ \mathrm{Gr}_0^P H^{j}(h^{-1}(U)) &&\quad\\
&=\dim_\bQ H^j(U)&&\quad(\text{by Remark~\ref{rmk:fullsupp}})\\
&=\dim_\bQ H^j(|\beta|)&&\quad(\text{by \cite[Lemma~2.4]{PSSZ24}})\\
&=\begin{cases}
    1&\text{(if $j$ is even)}\\
    0&\text{(if $j$ is odd)}
\end{cases},&&\quad
\end{aligned}\]
which proves the induction base. Then suppose that $n_{h_U}^{i',j'}=n_{\infty}^{i',j'}$ holds for any $i',j'\in\bZ_{\ge0}$ satisfying $i'< i$ and $i'+j'\le i+j$. For any $\ell\le i+j$, it follows from \eqref{eq:myms}, our induction hypothesis, and \eqref{eq:sk} that
\[\begin{aligned}n_{h_U}^{i,\ell-i}&=\sum_{i'+j'=i;\,i',j'\ge0}n_{h_U}^{i',\ell-i'-2j'}-\sum_{i'<i;\,i'+j'=i;\,i',j'\ge0}n_{h_U}^{i',\ell-i'-2j'}\\
&=\left(b_\ell(\cC^{[i]}_{U})-b_{\ell}(\cC^{[i-1]}_{U})\right)-\sum_{i'<i;\,i'+j'=i;\,i',j'\ge0}n_\infty^{i',\ell-i'-2j'}\\
&=\sum_{r=0}^{\lfloor \ell/2\rfloor}\left(b_{\ell-2r}(S^{[i]})-b_{\ell-2r}(S^{[i-1]})\right)-\sum_{i'<i;\,i'+j'=i;\,i',j'\ge0}n_\infty^{i',\ell-i'-2j'}.
\end{aligned}\]
To complete the induction step, it suffices to show
\begin{equation}\label{eq:suff}
\sum_{r=0}^{\lfloor \ell/2\rfloor}\left(b_{\ell-2r}(S^{[i]})-b_{\ell-2r}(S^{[i-1]})\right)=\sum_{k=0}^i n_\infty^{k,\ell-2i+k},
\end{equation}
which is equivalent to
\[\text{$z^\ell w^i$-coefficient of $G(z,w)\cdot\frac{1-w}{1-z^2}$}=\text{$q^i t^{\ell-i}$-coefficient of $\frac{H(q,t)}{1-qt}$}\]
with the notation as in Remark~\ref{rmk}.
Hence \eqref{eq:suff} follows from \eqref{eq:GH}.
\end{proof}

Finally, we prove our first main theorem stated in Section~\ref{sec:intro} by taking $U$ to be $U_k$ as given in Theorem~\ref{thm:gsm}.

\begin{proof}[Proof of Theorem~\ref{thm:main}]
By the upper semicontinuity of fiber dimension, we have
\begin{equation}
\label{eq:c}
    N(\beta)\ge 2c-2.
\end{equation}
    Since $\beta$ is $\max\{1,2k-2\}$-very ample, let $U_k\subset |\beta|^{\mathrm{int}}$ be as in Theorem~\ref{thm:gsm}. Using \eqref{eq:c}, the assumption $k\le c-1$, and Theorem~\ref{thm:gsm}, we obtain the following properties:
    \begin{enumerate}
    \item $N(U_k)\ge\min\{N(\beta),2k\}=2k$;
    \item $h^{-1}(U_k)$ is smooth.
\end{enumerate}
Note that $|\beta|$ contains a connected smooth curve by Bertini's theorem. Applying Corollary~\ref{cor} with $U=U_k$, we have
\begin{equation}
\label{eq:bk}
b_k(h^{-1}(U_k))=b_k^\infty.
\end{equation}
Since fibers of $h$ over integral curves have the same dimension, Theorem~\ref{thm:gsm} implies 
\[\codim(h^{-1}(|\beta|^{\mathrm{int}})\setminus h^{-1}(U_k),h^{-1}(|\beta|^{\mathrm{int}}))=\codim(|\beta|^{\mathrm{int}}\setminus U_k,|\beta|^{\mathrm{int}})\ge k+1.\]
This, together with the assumption on $c=\codim(M_{\beta,\chi} \setminus h^{-1}(|\beta|^{\mathrm{int}}),M_{\beta,\chi})$, yields
\[\codim(M_{\beta,\chi}\setminus h^{-1}(U_k),M_{\beta,\chi})\ge \min\{c,k+1\} > k.\]
Then by \cite[Theorem~6.7.4]{Max19} and the smoothness of $h^{-1}(U_k)$, we have isomorphisms
\[\IH^k(M_{\beta,\chi})\cong\IH^k(h^{-1}(U_k))\cong H^k(h^{-1}(U_k)),\]
and therefore, by \eqref{eq:bk},
\[\dim_\bQ \IH^k(M_{\beta,\chi})=b_k(h^{-1}(U_k))=b_{k}^\infty,\]
which completes the proof.
\end{proof}

\subsection{Discussion on the general case}
When $\beta$ is indecomposable (i.e., not linearly equivalent to $C_1+C_2$ for any nonzero effective divisors $C_1$ and $C_2$), all curves in $|\beta|$ are integral and $M_{\beta,\chi}$ is irreducible for any $\chi\in\bZ$. Now we provide an example, with $S$ a surface of general type and $\beta$ indecomposable, where Theorem~\ref{thm:main} applies to large $k$.

\begin{example}\label{ex:gentype}
For $m\ge 2$, write
\[N\vcentcolon=\dim (\bP^{m-1})^{[m+1]}+2m+3,\]
where $(\bP^{m-1})^{[m+1]}$ is the Hilbert scheme of $m+1$ points on $\bP^{m-1}$.
Let $S$ be a very general smooth complete intersection of $N-2$ hypersurfaces of degree $d\ge m$ in $\bP^N$, and let $\beta\in|\cO_{\bP^N}(1)|_S|$. Then $\beta$ is indecomposable since $\Pic(S)=\bZ[\cO_S(\beta)]$ by \cite[Theorem 1]{Kim91}. We claim that $\beta$ is $m$-very ample. 
If $Z\subset \bP^N$ is a $0$-dimensional closed subscheme of length $m+1$ such that the restriction map
\begin{equation}\label{eq:nsurj}
r_Z':H^0(\bP^N,\cO_{\bP^N}(1))\to H^0(Z,\cO_{\bP^N}(1)|_Z)
\end{equation}
is not surjective, then $Z$ is contained in some linear subspace $\L\subset\bP^N$ of dimension $m-1$. 
Hence, the scheme $B$ parametrizing those $Z$ satisfies
\[\dim B\le \dim \operatorname{Gr}(m, N+1)+\dim(\bP^{m-1})^{[m+1]}=m(N+1-m)+\dim(\bP^{m-1})^{[m+1]}.\]
By \cite[Example~1.8.30]{Laz04}, the restriction map
\[\tilde r_Z: H^0(\bP^N,\cO_{\bP^N}(d))\to H^0(Z,\cO_{\bP^N}(d)|_Z)\]
is surjective, and thus $\dim(\ker\tilde r_Z)=h^0(\bP^N,\cO_{\bP^N}(d))-(m+1)$. The incidence correspondence
\[\begin{tikzcd}
	{\{(H_1,\cdots,H_{N-2},Z)\in |\cO_{\bP^N}(d)|^{N-2}\times B:Z\subset H_1\cap\cdots\cap H_{N-2}\}} & B \\
	{|\cO_{\bP^N}(d)|^{N-2}}
	\arrow["p_1", from=1-1, to=1-2]
	\arrow["p_2"', from=1-1, to=2-1]
\end{tikzcd}\]
together with the inequality
\[\dim B\le m(N+1-m)+\dim(\bP^{m-1})^{[m+1]}<(N-2)(m+1)\]
implies that $p_2$ is not surjective, so we can arrange that $S$ does not contain any $Z\in B$. Since $H^0(S,\cO_S(\beta))=H^0(\bP^N,\cO_{\bP^N}(1))$, we conclude that $\beta$ is $m$-very ample. It follows from Theorem~\ref{thm:main} that 
\[\dim_\bC\IH^k(M_{\beta,\chi})=b_k^\infty\]
holds for $k\le m/2+1$.
\end{example}

In general, the moduli space $M_{\beta,\chi}$ can be reducible even if $\beta$ is sufficiently positive (see \cite[\S 5]{SZ25}). Moreover, even if we assume that $M_{\beta,\chi}$ is irreducible, it might not be true that the codimension $c$ in Theorem~\ref{thm:main} is the same as $\codim(|\beta|\setminus|\beta|^{\mathrm{int}},|\beta|)$. 

In the cases where $M_{\beta,\chi}$ is irreducible and $c=\codim(|\beta|\setminus|\beta|^{\mathrm{int}},|\beta|)$, the asymptotic behavior of $\IH^k(M_{\beta,\chi})$ can be described more precisely, provided that $\codim(|\beta|\setminus|\beta|^{\mathrm{int}},|\beta|)$, or equivalently $N(\beta)$ (defined in \S \ref{sec:step}), can be calculated.
The following lemma allows us to estimate $N(\beta)$ by analyzing the decompositions of $\beta$.

\begin{lem}
\label{lem:fin}
Fix an ample divisor $\beta$ on $S$. There are at most finitely many pairs of classes $(\theta_1,\theta_2)$ in $\Num(S)$ such that $\theta_i$ ($i=1,2$) are the classes of nonzero effective divisors $C_i$ with
\[C_1+C_2=\beta.\]
\end{lem}
\begin{proof}
Let $D_1,\cdots,D_{\rho}$ be an orthogonal basis for $\Num(S)_\bQ\vcentcolon= \Num(S)\otimes_\bZ \bQ$, where $D_1$ is the class of an ample divisor. By the Hodge index theorem, $D_\ell^2<0$ for $1<\ell\le \rho$.
Given any $C_1$, $C_2$ as in the statement, write 
\[\beta\equiv\sum_{j=1}^{\rho}a_jD_j,\quad C_i\equiv\sum_{j=1}^{\rho}a_{i,j}D_j\]
for some $a_j,a_{i,j}\in\bQ$, where $\equiv$ denotes the numerical equivalence. For $1<\ell\le \rho$, let $n_\ell>0$ be an integer such that $A_\ell\vcentcolon=n_\ell D_1+D_\ell$ is ample. Then $D_1\cdot C_i>0$ and $A_\ell\cdot C_i>0$, which implies
\[a_{i,1}>0\quad\text{and}\quad n_\ell a_{i,1}D_1^2+a_{i,\ell}D_\ell^2>0 .\]
Combined with $D_\ell^2<0$ and $C_1+C_2=\beta$, this shows that the possible values of $a_{i,j}$ ($i=1,2$; $j=1,\cdots,\rho$) are bounded since $\beta$ is fixed, and thus the result follows.
\end{proof}

The next proposition shows that $N(\beta)$ can be arbitrarily large as $\beta$ becomes positive enough, generalizing \cite[Proposition~1.5]{PSSZ24}.

\begin{prop}\label{prop:Nbeta}
For any ample divisor $\beta$ on $S$, 
\[\lim_{d\to+\infty}N(d\beta)=+\infty.\]
\end{prop}
\begin{proof}
Any non-integral curve in $|d\beta|$ is in the image of a morphism from a scheme $P$ to $|d\beta|$, where $P$ parametrizes the triples $(C,D,N)$ satisfying $N\in\Pic^0(S)$, $C\in|C_1+N|$, and $D\in|C_2-N|$ for two nonzero effective divisors  $C_1$, $C_2$ such that $C_1+C_2=d\beta$. By Lemma~\ref{lem:fin}, the images of finitely many such morphisms cover $|d\beta|\setminus|d\beta|^{\mathrm{int}}$. Thus,
\[N(d\beta)\ge2\min_{C_1+C_2=d\beta}\left\{\dim|d\beta|-\dim |C_1|\times |C_2|\times\Pic^0(S)\right\}-2,\]
so it suffices to estimate all possible values of
\[\dim|d\beta|-\dim |C_1|\times |C_2|\times\Pic^0(S)=\dim|d\beta|-\dim|C_1|-\dim|C_2|-h^1(\cO_S).\]

Choose $d_0\in\bZ_{>0}$ such that $H_0\vcentcolon=d_0\beta$ is very ample, and fix a smooth curve $C_0\in |H_0|$.
We claim that there exists a constant $A>0$, depending only on $S$ and $H_0$, such that for every nonzero effective divisor $D$, 
\begin{equation}\label{eq:claim}
\dim|D|\le \frac{(D\cdot H_0)^2}{2H_0^2}-\frac{(K_S\cdot H_0)(D\cdot H_0)}{2H_0}+A.
\end{equation}
To prove the claim,  write 
\[\ell\vcentcolon=K_S\cdot H_0,\quad s\vcentcolon=H_0^2, \quad t\vcentcolon=D\cdot H_0=as+r\]
with $a,r\in\bZ_{\ge0}$ and $0\le r<s$. Since $(D-(a+1)H_0)\cdot H_0<0$, the line bundle $\cO_S(D-(a+1)H_0)$ has no nonzero sections. Repeatedly using the exact sequence
\[0\to \cO_S(D-(i+1)H_0)\to \cO_S(D-iH_0)\to \cO_S(D-iH_0)|_{C_0}\to 0\]
yields
\begin{equation}\label{eq:dimD}
h^0(S,\cO_S(D))\le \sum_{i=0}^{a} h^0\left(C_0, \cO_S(D-iH_0)|_{C_0}\right).
\end{equation}
Since the degree of $\cO_S(D-iH_0)|_{C_0}$ is $t-is$ and the genus of $C_0$ is $1+(s+\ell)/2$, it follows from the Riemann–Roch formula and \cite[Ch.~IV, Example 1.3.4]{Har77} that if $t-is>2+s+\ell$, then
\begin{equation}\label{eq:lindim1}
h^0\left(C_0, \cO_S(D-iH_0)|_{C_0}\right)=\chi(\cO_S(D-iH_0)|_{C_0})=t-is-\frac{s+\ell}{2}.
\end{equation}
By Clifford’s theorem \cite[Ch.~IV, Theorem 5.4]{Har77}, if $t-is\le 2+s+\ell$, then
\begin{equation}\label{eq:lindim2}
h^0\left(C_0, \cO_S(D-iH_0)|_{C_0}\right)\le \frac{t-is}{2}+1.
\end{equation}
The number of elements in $\{i\in\bZ_{\ge0}:t-is\le 2+s+\ell\}$ is at most $n_0$, where $n_0$ depends only on $S$ and $H_0$. Thus, by \eqref{eq:lindim1} and \eqref{eq:lindim2}, we have
\[\begin{aligned}
\sum_{i=0}^{a} h^0\left(C_0, \cO_S(D-iH_0)|_{C_0}\right)&\le\sum_{i=0}^{a}\left(t-is-\frac{s+\ell}{2}\right)+\left(1+\frac{s+\ell}{2}\right)n_0\\
&=\frac{t^2}{2s}-\frac{\ell t}{2s}+\left(r-\frac{r^2}{2s}+\frac{r\ell }{2s}-\frac{s+\ell }{2}\right)+\left(1+\frac{s+\ell }{2}\right)n_0.
\end{aligned}\]
Combining this with \eqref{eq:dimD} proves \eqref{eq:claim}.

By the Riemann–Roch formula and Serre vanishing, for $d\gg 0$,
\[
\dim|d\beta|= \chi(\cO_S(d\beta))-1=\frac{(d\beta\cdot H_0)^2}{2H_0^2}-\frac{(K_S\cdot H_0)(d\beta\cdot H_0)}{2H_0^2}+\chi(\cO_S)-1.
\]
Therefore, applying \eqref{eq:claim} to $C_1$ and $C_2$, we conclude that for $d\gg0$,
\[\begin{aligned}
&\hphantom{{}={}}\!\dim|d\beta|-\dim|C_1|-\dim|C_2|-h^1(\cO_S)\\
&\ge\frac{(d\beta\cdot H_0)^2}{2H_0^2}-\frac{(K_S\cdot H_0)(d\beta\cdot H_0)}{2H_0^2}+\chi(\cO_S)-1-h^1(\cO_S)\\
&\hphantom{{}={}}-\left[\frac{(C_1\cdot H_0)^2}{2H_0^2}-\frac{(K_S\cdot H_0)(C_1\cdot H_0)}{2H_0^2}+A+\frac{(C_2\cdot H_0)^2}{2H_0^2}-\frac{(K_S\cdot H_0)(C_2\cdot H_0)}{2H_0^2}+A\right]\\
&=\frac{(C_1\cdot H_0)(d\beta\cdot H_0-C_1\cdot H_0)}{H_0^2}+\chi(\cO_S)-1-h^1(\cO_S)-2A\\
&\ge\frac{d\beta\cdot H_0-1}{H_0^2}+\chi(\cO_S)-1-h^1(\cO_S)-2A,
\end{aligned}\]
which completes the proof.
\end{proof}

\section{Stabilization for surfaces with nef anticanonical divisor}
\label{sec:stab}
Throughout this section, except in \S \ref{subsec:ques}, we assume that $-K_S$ is nef. 

\subsection{Proof of stabilization}
\subsubsection{Stabilization of intersection Betti numbers}
Together with Proposition~\ref{prop:irr} and Proposition~\ref{prop:Nbeta}, Theorem~\ref{thm:main} implies the stabilization of the intersection Betti numbers of $M_{\beta,\chi}$.
\begin{proof}[Proof of Theorem~\ref{thm:stab}]
Choose $d(\beta_0,k)\in\bZ_{>0}$ such that for all $d\ge d(\beta_0,k)$, the following properties hold:
\begin{enumerate}
    \item $N(d\beta_0)\ge 2k$ (guaranteed by Proposition~\ref{prop:Nbeta});
    \item $d\beta_0$ is $\max\{1,2k-2\}$-very ample;
    \item $2\dim|d\beta_0|\ge 3k$.
\end{enumerate} 
Since $-K_S$ is nef, the divisor $d\beta_0-K_S$ is ample and thus $\dim|d\beta_0|=\chi(\cO_S(d\beta_0))-1$ by the Kodaira vanishing theorem.
By Proposition~\ref{prop:irr} (i), $M_{d\beta_0,\chi}$ is irreducible, and by Proposition~\ref{prop:irr} (ii), we have
\begin{equation*}
\label{eq:codimeq}
\begin{aligned}
   \codim(M_{d\beta_0,\chi}\setminus h^{-1}(|d\beta_0|^{\mathrm{int}}),M_{d\beta_0,\chi})&=\codim(|d\beta_0|\setminus |d\beta_0|^{\mathrm{int}},|d\beta_0|)\\
   &=1+\frac{1}{2}N(d\beta_0)\ge k+1. 
\end{aligned}
\end{equation*}
Hence, the result follows from  Theorem~\ref{thm:main} with $\beta=d\beta_0$.
\end{proof}



\subsubsection{Stabilization of perverse Hodge numbers}
Suppose that $\beta$ is a base-point-free, ample divisor on $S$ and $M_{\beta,\chi}$ is smooth. Then $|\beta|$ contains a connected smooth curve by Bertini's theorem \cite[Ch.~III, Corollary~10.9, Remark~10.9.1]{Har77}. By the Kodaira vanishing theorem, Proposition~\ref{prop:irr} applies to this situation. 
Recall that $n_{\beta,\chi}^{i,j}$ is the number defined by \eqref{eq:nij} with $\varphi$ taken to be $h:M_{\beta,\chi}\to|\beta|$. 

In order to apply Theorem~\ref{thm} to compute $n_{\beta,\chi}^{i,j}$, we need the following lemma which is parallel to \cite[Corollary~5.3]{MSY25}. For the reader’s convenience, we provide a proof below.

\begin{lem}
\label{lem:intEnr}
Let $\beta$ be a base-point-free, ample divisor on a surface $S$ with $-K_S$ nef. If $M_{\beta,\chi}$ is smooth, then for any $i,j\in\bZ_{\ge0}$ satisfying $i+j\leq N(\beta)$, 
\[n_{\beta,\chi}^{i,j}=n^{i,j}_{h^{\circ}},\]
where $h^\circ=h_{|\beta|^{\mathrm{int}}}:h^{-1}(|\beta|^{\mathrm{int}})\to|\beta|^{\mathrm{int}}$ is the restriction of $h:M_{\beta,\chi}\to|\beta|$.
\end{lem}
\begin{proof}
    Let $W=M_{\beta,\chi}\setminus h^{-1}(|\beta|^{\mathrm{int}})$. Denote by $\iota:W\hookrightarrow M_{\beta,\chi}$ the closed embedding, and by $\jmath:h^{-1}(|\beta|^{\mathrm{int}}) \hookrightarrow M_{\beta,\chi}$, $\jmath^{\circ}:|\beta|^{\mathrm{int}}\hookrightarrow|\beta|$ the open embeddings. By Proposition~\ref{prop:irr} (ii), 
$$\codim(W,M_{\beta,\chi})=\codim(|\beta|\setminus|\beta|^{\mathrm{int}},|\beta|)=1+\frac{1}{2}N(\beta).$$
Consider the exact triangle (all functors are derived)
\[
\iota_*\iota^!\bQ_{M_{\beta,\chi}} \longrightarrow \bQ_{M_{\beta,\chi}} \longrightarrow \jmath_{*}\jmath^{*}\bQ_{M_{\beta,\chi}} \stackrel{+1}\longrightarrow.
\]
By the property of Verdier duals \cite[Proposition~5.3.9]{Max19} and the smoothness of $M_{\beta,\chi}$, 
\[\iota_*\iota^!\bQ_{M_{\beta,\chi}}\simeq (D_{M_{\beta,\chi}}\iota_*\bQ_W)[-2\dim M_{\beta,\chi}],\]
and therefore for $k\le N(\beta)+1$, 
\[\bH^{k}(M_{\beta,\chi},\iota_*\iota^!\bQ_{M_{\beta,\chi}})\cong H^{2\dim M_{\beta,\chi}-k}(W)^\vee=0.\]
By the long exact sequence in hypercohomology, the vanishing of $\bH^{k}(M_{\beta,\chi},\iota_*\iota^!\bQ_{M_{\beta,\chi}})$ shows that for $\ell\le N(\beta)$, the restriction map induces an isomorphism
 \begin{equation}
\label{eq:cohiso}
    H^\ell(M_{\beta,\chi})\cong H^\ell(h^{-1}(|\beta|^{\mathrm{int}})).
\end{equation}
Since the restriction to $h^{-1}(|\beta|^{\mathrm{int}})$ clearly preserves a relatively ample class (relative to $h$), it follows from \cite[Lemma~3.3]{dCMS22} that $n_{\beta,\chi}^{i,j}\le n_{h^\circ}^{i,j}$ for $i+j\le N(\beta)$. Then the result follows by \eqref{eq:nijsum} and \eqref{eq:cohiso}.
\end{proof}

Now we can prove the stabilization of $n_{\beta,\chi}^{i,j}$.


\begin{proof}[Proof of Theorem~\ref{thm:Enrs}]
Choose $d(\beta_0,i,j)$ such that for all $d\ge d(\beta_0,i,j)$, the following properties hold:
\begin{enumerate}
    \item $N(d\beta_0)\ge i+j$ (guaranteed by Proposition~\ref{prop:Nbeta});
    \item $d\beta_0$ is $\max\{0,i-1\}$-very ample;
    \item $2\dim|d\beta_0|=2\chi(\cO_S(d\beta_0))-2\ge 3i+j$.
\end{enumerate}
Applying Lemma~\ref{lem:intEnr} and Theorem~\ref{thm} with $\beta=d\beta_0$ and $U=|d\beta_0|^{\mathrm{int}}$, we have
$$n_{d\beta_0,\chi}^{i,j}=n_{h_U}^{i,j}=n_\infty^{i,j},$$
which completes the proof.
\end{proof}

\subsubsection{The case of Enriques surfaces}
Motivated by \cite{Obe24}, we explore the case of Enriques surfaces in more detail.

Recall that an Enriques surface $S$ is a minimal surface of Kodaira dimension $0$ with $h^1(\cO_S)=h^2(\cO_S)=0$. Its second Betti number is $b_2(S)=10$. Assume that $S$ is a generic Enriques surface. For an ample divisor $\beta$ whose class in $\Num(S)$ is not divisible by $2$, the moduli space $M_{\beta,1}$ is a smooth projective Calabi--Yau variety by \cite{Sac19}. 
The linear system $|\beta|$ contains a connected smooth curve (e.g., \cite[Corollary~1.13]{Sac19}).

%

The following property of linear systems on $S$ is well-known (e.g., \cite[Theorem]{KL07}) and was used in \cite[Proposition~1.14]{Sac19} to compute the codimension of the locus of reducible curves in a linear system.
\begin{lem}\label{lem:Enr}
Let $D$ be a nonzero effective nef divisor  on an Enriques surface $S$. 
    \begin{enumerate}[(i)]
        \item If $D^2>0$, then \[\dim|D|=\chi(\cO_S(D))-1=\frac{D^2}{2}.\]
    \item If $D^2=0$, then $D$ is linearly equivalent to either $mE$ or $mE+K_S$ for some $m\in\bZ_{>0}$ and some primitive effective divisor $E$ with $E^2=0$. Moreover,
    \[\dim|D|=\begin{cases}
        \lfloor\frac{m}{2}\rfloor&(\text{if }D\sim mE)\\
        \lfloor\frac{m-1}{2}\rfloor&(\text{if }D\sim mE+K_S).
    \end{cases}\]
    \end{enumerate}
\end{lem}

To obtain an explicit expression for $d(\beta_0,i,j)$ in the case of generic Enriques surfaces, we use Lemma~\ref{lem:Enr} to bound $N(d\beta_0)$. By the proof of Proposition~\ref{prop:Nbeta}, it suffices to estimate all possible values of
\[\dim|d\beta_0|-\dim |C_1|\times |C_2|\times\Pic^0(S)=\dim|d\beta_0|-\dim|C_1|-\dim|C_2|.\]
Note that $C_1$ and $C_2$ are nef since $S$ is a generic Enriques surface.

\noindent\textit{\emph{Case 1.1}: $C_1^2>0$ and $C_2^2>0$.} 
    In this case, $C_1^2\ge2$ and $C_2^2\ge2$. By Lemma~\ref{lem:Enr} (i) and the Hodge index theorem, we have
    \[\begin{aligned}
        \dim|d\beta_0|-\dim|C_1|-\dim|C_2|&=C_1\cdot C_2\ge \sqrt{C_1^2C_2^2}\\
        &\ge \sqrt{2(d^2\beta_0^2-2C_1\cdot C_2-2)},
    \end{aligned}   
    \]
    where the last inequality follows from $d^2\beta_0^2=C_1^2+C_2^2+2C_1\cdot C_2$. Therefore,
    \begin{equation}
    \label{eq:1.1}
    \dim|d\beta_0|-\dim|C_1|-\dim|C_2|=C_1\cdot C_2\ge d\sqrt{2\beta_0^2}-2.
    \end{equation}

\noindent\textit{\emph{Case 1.2}: $C_1^2=0$ and $C_2^2>0$.} By Lemma~\ref{lem:Enr} (ii), $C_1$ is linearly equivalent to $m_1E_1$ or $m_1E_1+K_S$ for some $m_1\in\bZ_{>0}$ and some primitive effective divisor $E_1$ with $E_1^2=0$. Moreover,
    \begin{equation}
    \label{eq:1.2}
        \begin{aligned}
      \dim|d\beta_0|-\dim|C_1|-\dim|C_2|&\ge C_1\cdot C_2-\lfloor \frac{m_1}{2}\rfloor\\ 
       &\ge m_1\left(d\beta_0\cdot E_1-\frac{1}{2}\right)\ge d-\frac{1}{2}.
    \end{aligned}
    \end{equation}

\noindent\textit{\emph{Case 1.3}: $C_1^2=C_2^2=0$.} Let $m_1$ and $E_1$ be as in Case 1.2. By Lemma~\ref{lem:Enr} (ii), $C_2$ is linearly equivalent to $m_2E_2$ or $m_2E_2+K_S$ for some $m_2\in\bZ_{>0}$ and some primitive effective divisor $E_2$ with $E_2^2=0$. Then we have $C_1\cdot C_2=m_1m_2E_1\cdot E_2=d^2\beta_0^2/2>0$ and
\begin{equation}
\label{eq:1.3}
    \begin{aligned}
      \dim|d\beta_0|-\dim|C_1|-\dim|C_2|&\ge C_1\cdot C_2-\lfloor \frac{m_1}{2}\rfloor-\lfloor \frac{m_2}{2}\rfloor\ge \frac{d^2\beta_0^2-(m_1+m_2)}{2}\\
      &\ge \frac{d^2\beta_0^2-(1+d^2\beta_0^2/2)}{2}=\frac{d^2\beta_0^2-2}{4}.
    \end{aligned}
\end{equation}

By \eqref{eq:1.1}--\eqref{eq:1.3}, \cite[Proposition~2.5]{Sze01}, and the proof of Theorem~\ref{thm:Enrs}, we can take
\begin{equation}
\label{eq:d0}
    d(\beta_0,i,j)=\max\left\{i+1,\,\left\lceil\frac{i+j+2}{2}\right\rceil,\,\left\lceil\frac{i+j+6}{2\sqrt{2\beta_0^2}}\right\rceil,\,\left\lceil\sqrt{\frac{2i+2j+6}{\beta_0^2}}\,\right\rceil\right\}
\end{equation}
in Theorem~\ref{thm:Enrs}. 


\subsection{Some related stabilization questions}
\label{subsec:ques}
We collect some stabilization questions for moduli spaces of one-dimensional sheaves. Now $S$ is not necessarily a surface whose anticanonical divisor is nef, provided the conjecture and questions below make sense. We begin by providing a refined version of Conjecture~\ref{conj}, stated heuristically as follows.
\begin{conj}\label{conj:stabPHN}
  For fixed $i,j\in\bZ_{\ge0}$ and $\chi\in\bZ$, the perverse Hodge number $n_{\beta,\chi}^{i,j}$ stabilizes as $\beta$ becomes sufficiently positive. Furthermore, $n_{\beta,\chi}^{i,j}$ stabilizes to $n_\infty^{i,j}$ if $M_{\beta,\chi}$ is asymptotically irreducible.
\end{conj}

In \cite{PSSZ24}, the approach to Conjecture~\ref{conj} for del Pezzo surfaces is to prove Conjecture~\ref{conj:stabPHN}. When $S$ is a K3 or an Abelian surface and $M_{\beta,\chi}$ is smooth, Shen--Yin \cite{SY22} proved that the perverse Hodge numbers of $M_{\beta,\chi}$ are equal to its Hodge numbers using hyper-Kähler geometry and thus the stabilization of the perverse Hodge numbers follows from the closed formula of Hodge numbers of $S^{[n]}$ \cite[Theorem~2.3.14]{Got06}.  Our Theorem~\ref{thm:Enrs} also shows the stabilization of perverse Hodge numbers when $-K_S$ is nef and $M_{\beta,\chi}$ is smooth.

\begin{question}\label{IHN}
For fixed $i,j\in\bZ_{\ge0}$ and $\chi\in\bZ$, does the intersection Hodge number $\dim_\bQ \IH^{i,j}(M_{\beta,\chi})$ stabilize as $\beta$ becomes sufficiently positive?
\end{question}
In \cite{MM24}, Mauri and Migliorini showed that Question~\ref{IHN} has a positive answer for K3 surfaces, based on the previously known case where $\chi=1$ and $S$ is a K3 surface. 
By \cite[Theorem~0.2]{Bou22}, in the case $S=\bP^2$, Question~\ref{IHN} is reduced to the question of stabilization of intersection Betti numbers, and thus has a positive answer by \cite{Yua23} or \cite{PSSZ24}. It seems that our method in this paper does not work for intersection Hodge numbers in the new cases we prove for Conjecture~\ref{conj}.

The stabilization problems for moduli spaces of one-dimensional sheaves on $S$ can also be studied from a motivic perspective. Let $K_0(\mathrm{Var})$ be the Grothendieck ring of varieties, which is the free Abelian group generated by (the isomorphism classes of) varieties $[X]$ modulo the scissor relations
\[[X]=[Y]+[Z],\]
where $Y$ and $Z$ are disjoint locally closed subvarieties of $X$ with $X=Y\cup Z$. The product in $K_0(\mathrm{Var})$ is defined by taking the fiber product.
Let $\mathbb{L}\in K_0(\mathrm{Var})$ be the class of the affine line $\bA^1$.
Denote by $\widehat{K}_0(\mathrm{Var})$ the completion of $K_0(\mathrm{Var})[\mathbb{L}^{-1}]$ with respect to the dimension filtration (see \cite[p.~17]{CW22} for more details). 
\begin{question}\label{ques:motivic}
    For fixed $\chi\in\bZ$, does the class $[M_{\beta,\chi}]$ stabilize in $\widehat{K}_0(\mathrm{Var})$ as $\beta$ becomes sufficiently positive?
\end{question}

Coskun and Woolf \cite{CW22} proved that when $S$ is a rational surface and $K_S\cdot H<0$, the class of the moduli space of $H$-stable torsion-free sheaves with fixed Chern character on $S$ stabilizes in $\widehat{K}_0(\mathrm{Var})$ as the determinant of sheaves tends to $+\infty$. A positive answer to Question \ref{ques:motivic} would establish an analogous result in the rank $0$ case.

A natural further question is whether there is a geometric interpretation for the stabilization of Betti numbers. One potential approach is to consider the tautological classes in $H^*(M_{\beta,\chi})$, i.e., the Künneth components of Chern classes of the universal sheaf (assuming its existence) on $M_{\beta,\chi}\times S$. Let $RH^*(M_{\beta,\chi})\subset H^*(M_{\beta,\chi})$ be the subring generated by the tautological classes. Although $RH^*(M_{\beta,\chi})\neq H^*(M_{\beta,\chi})$ in general (cf.~\cite[Theorem~5.10]{Sac19}), we expect the stabilization of Betti numbers of $M_{\beta,\chi}$ to be explained by the stabilization of $RH^*(M_{\beta,\chi})$. When $S$ is a del Pezzo surface and $M_{\beta,\chi}$ is smooth, one has $RH^*(M_{\beta,\chi})=H^*(M_{\beta,\chi})$ and the stable Betti numbers are determined by the combinatorics of tautological classes (see \cite[Corollary 1.8]{PSSZ24}). If our expectation is true, the stabilization phenomenon in the del Pezzo case fits into a more general picture.


As for the intersection cohomology, the $k$-th stable intersection Betti number calculated in this paper also coincides with $b_k^\infty$. 
In what follows, we assume that $M_{\beta,\chi}$ is pure dimensional.
Note that there is a homomorphism of graded vector spaces (e.g., \cite[Exercise~6.7.1, Remark~6.7.2]{Max19})
\[ \alpha:  H^*(M_{\beta,\chi}) \rightarrow  \IH^*(M_{\beta,\chi}). \] Thus we can  define $ R\IH^\ast(M_{\beta,\chi})$ as the image of $RH^\ast (M_{\beta,\chi})$ under the map $\alpha$. Given our results, we also expect the stabilization of intersection Betti numbers to be explained by the stabilization of $ {R\IH}^\ast(M_{\beta,\chi})$, and therefore we ask the following question.  
\begin{question} Is it true that $\dim_\bQ {\RIH}^k(M_{\beta,\chi})=b_k^\infty$ if $\beta$ is sufficiently positive?
\end{question}

Another direction is to explore the tautological relations in $RH^\ast(M_{\beta,\chi})$. When $S=\bP^2$, the tautological relations are studied in \cite{KLMP24}, leading to the formulation of a strengthened version of the $P=C$ conjecture proposed in \cite{KPS23}. For other surfaces, however, the problem remains unexplored.

More generally, one may consider similar questions for moduli spaces of one-dimensional sheaves on a projective surface with mild singularities. In this setting, we continue to use the notation $M_{\beta,\chi}$. 

\begin{question}
\label{que:sing}
Fix $\chi\in\bZ$ and let $\beta$ be a divisor on a simply connected projective surface with ADE singularities. For fixed $k \in\bZ_{\ge 0}$, does the intersection Betti number $\dim_\bQ \IH^{k}(M_{\beta,\chi})$ stabilize as $\beta$ becomes sufficiently positive? Similarly, for fixed $i,j \in\bZ_{\ge 0}$, does the intersection Hodge number $\dim_\bQ \IH^{i,j}(M_{\beta,\chi})$ stabilize?
\end{question}

Toward answering Question~\ref{que:sing}, we expect our strategy in this paper to work in some cases after suitable modifications, but there are some difficulties. For example, the relative Hilbert schemes of points over the linear system $|\beta|$ may not be equidimensional as the curves in $|\beta|$ may not be locally planar. Moreover,  the Betti numbers (or intersection Betti numbers) of Hilbert schemes of points on a singular surface are still unknown.

\section*{Acknowledgments} 
The authors would like to thank Mirko Mauri, Weite Pi, and Junliang Shen for helpful discussions. FZ is grateful to her advisor Jun Li for his constant support, to János Kollár and Claire Voisin for inspiring discussions on the deformation theory of curve singularities, to Georg Oberdieck for suggesting considering the case of bielliptic surfaces, and to Songsong Huang for useful conversations. 

FS is partially supported by the Fundamental Research Funds for the Central Universities, NSFC (No. 12025106) and Shaanxi NSF (No. 2025JC-QYCX-002). FZ is supported by the NSFC grants (No. 12121001 and No. 12425105) and LMNS.

\end{document}